\documentclass[conference]{IEEEtran}
\IEEEoverridecommandlockouts
\usepackage[dvipsnames]{xcolor}
\usepackage{cite}
\usepackage{amsmath,amssymb,amsfonts,amsthm}
\usepackage{algorithm}
\usepackage{algorithmic}
\usepackage{graphicx}
\usepackage{textcomp}
\usepackage[hidelinks]{hyperref}
\usepackage{mleftright}
\usepackage{physics}
\usepackage{pgfplots}
\usepackage{multirow}
\usepackage{dashbox}

\pgfplotscreateplotcyclelist{mycolor}{%
    red,mark=*\\%
    blue,mark=square\\%
    green,mark=star\\%
    orange,mark=triangle\\%
}
\pgfplotsset{
    compat=1.18,
    cycle list name=mycolor,
    table/col sep=comma,
    legend cell align=left,
    legend style={font=\scriptsize,draw=none,fill=none},
    every axis/.append style={font=\footnotesize}
}
\usepgfplotslibrary{groupplots}

\DeclareMathOperator{\lp}{fl} % Low precision
\DeclareMathOperator*{\argmin}{\arg\!\min}

\newcommand{\mchop}[1]{m_{\text{chop}_{#1}}}
\newcommand{\ump}{\widetilde{u}} % mixed precision u
\newcommand{\jac}{J}
\newcommand{\jacmp}{\widetilde{J}}
\newcommand{\fancyE}{\mathbf{\mathcal{E}}}
\newcommand{\KIOPS}{\textsc{KIOPS}}
\def\BibTeX{{\rm B\kern-.05em{\sc i\kern-.025em b}\kern-.08em
    T\kern-.1667em\lower.7ex\hbox{E}\kern-.125emX}}
\begin{document}

\title{Leveraging Mixed Precision in Exponential Time Integration Methods\\
\thanks{This work was performed under the auspices of the U.S. Department of Energy by Lawrence Livermore National Laboratory under Contract DE-AC52-07NA27344 and was supported by the LLNL-LDRD Program under Project No. 23-FS-013. LLNL-PROC-851497.}
}

\author{\IEEEauthorblockN{Cody J. Balos\IEEEauthorrefmark{1},
Steven Roberts\IEEEauthorrefmark{2}, and David J. Gardner\IEEEauthorrefmark{3}}
\IEEEauthorblockA{Center for Applied Scientific Computing,\\
Lawrence Livermore National Laboratory,\\
7000 East Ave, Livermore, CA\\ 
Email: \IEEEauthorrefmark{1}balos1@llnl.gov,
\IEEEauthorrefmark{2}roberts115@llnl.gov,
\IEEEauthorrefmark{3}gardner48@llnl.gov}}
\maketitle

% Temporarily turn page numbers on so we can see where we are at
\thispagestyle{plain}
\pagestyle{plain}

\begin{abstract}
The machine learning explosion has created a prominent trend in modern computer hardware towards low precision floating-point operations. 
In response, there have been growing efforts to use low and mixed precision in general scientific computing.
One important area that has received limited exploration is time-integration methods, which are used for solving differential equations that are ubiquitous in science and engineering applications.
In this work, we develop two new approaches for leveraging mixed precision in exponential time integration methods.
The first approach is based on a reformulation of the exponential Rosenbrock--Euler method allowing for low precision computations in matrix exponentials independent of the particular algorithm for matrix exponentiation.
The second approach is based on an inexact and incomplete Arnoldi procedure in Krylov approximation methods for
computing matrix exponentials and is agnostic to the chosen integration method.
We show that both approaches improve accuracy compared to using purely low precision and offer better efficiency than using only double precision when solving an advection-diffusion-reaction partial differential equation.
\end{abstract}

\begin{IEEEkeywords}
differential equations, mixed precision, high-performance computing
\end{IEEEkeywords}

\section{Introduction} 
In this paper we present two complementary concepts that enable accurate mixed precision computation in exponential time-integators. 
Exponential time-integrators are a class of numerical methods for solving ordinary differential equation (ODE) initial value problems of the form
\begin{equation} \label{eq:ODE}
    u'(t) = f(u(t)), \quad u(t_0) = u_0, \quad t \in [t_0, t_f],
\end{equation}
with $u(t) \in \mathbb{R}^N$. ODEs are ubiquitous across scientific domains and may arise directly from modeling some process or from discretizing a partial differential equation (PDE).
Exponential integrators are particularly well-suited to stiff problems due to their exact treatment of a linear term.
Alternative methods for stiff ODEs e.g., BDF or implicit Runge--Kutta methods, typically require an effective and efficient preconditioner which can be difficult to construct \cite{luan2017precon}.
Exponential time-integrators have been shown to be effective for many problems where practical preconditioners have not been developed \cite{hochbruck1998exponential, tokman2006efficient, hochbruck2009exponential,  loffeld2013compar, loffeld2014implementation, einkemmer2017performance}. 

Recent trends in computer hardware towards low precision floating point operations have been spurred largely by artificial intelligence and machine learning applications.  
Reuther et al. provides a comprehensive survey of current AI accelerators and their properties \cite{reuther2021ai}. 
The typical properties of this hardware indicate that leveraging low precision is necessary to achieve the full potential of much of this hardware.
% \note[id=DR]{The judgement call in the following sentence should be justified somehow; if this is a direct conclusion from [34] then consider combining the two sentences, e.g., “Reuther et al., … in [34], where it is clear that the use of … hardware.”}
This has resulted in significant interest in mixed precision computation.
The goal of incorporating mixed precision is to utilize the efficiency of low precision computation while maintaining an overall accuracy consistent with high precision computation.
What constitutes low and high precision depends on the context, but it is common to consider double precision as high precision and low precision as anything less than double. 
Mixed precision has been particularly popular in the numerical linear algebra \cite{abdelfattah2021survey,mukunoki2020performance,abdulah2021accelerating,haidar2020mixed,li2002design,baboulin2009accelerating,higham2022mixed,anzt2021towards} and deep learning literature \cite{uhlichmixed,micikeviciusmixed,nandakumar2020mixed,gong2019mixed,chen2021towards}.  
However, incorporating mixed precision into numerical time-integration methods has been studied much less \cite{grant2022perturbed,burnett2021performance,croci2022mixed}.
To better utilize current and emerging hardware capabilities, further research on incorporating mixed precision into numerical time-integration methods is needed. To this end, we present two approaches for leveraging mixed precision computations in exponential integrators:
\begin{enumerate}
    \item a reformulation of the exponential Rosenbrock--Euler method with order of accuracy $\order{h^2 + \epsilon h}$ instead of $\order{h^2 + \epsilon}$ where $\epsilon$ is the floating point precision,
    \item and the incorporation of low precision matrix-vector products in the Krylov approximation of matrix-exponential and vector products.
\end{enumerate}
These two approaches have different requirements and characteristics which may dictate which is most suitable for a particular application.
They can also be combined to create a practical and robust mixed precision exponential time-integrator.

The rest of this paper is organized as follows. In section \ref{s:expo_euler} we present the reformulated exponential Rosenbrock--Euler method.  This is followed by section \ref{s:KIIOPS} where we present the mixed precision Krylov approximation algorithm. 
In section \ref{s:tests} we demonstrate both approaches in solving an advection-diffusion-reaction PDE. 
Finally, in section \ref{s:conclusion} we provide key conclusions, impacts, and directions for future work. 
\section{Reformulation of Exponential Euler}
\label{s:expo_euler}

The exponential Rosenbrock--Euler method \cite{pope1963exponential} applied to \eqref{eq:ODE} is given by
\begin{equation} \label{eq:exp_Euler}
    u_{n+1} = u_n + h_n \varphi_1(h_n \jac_n) f(u_n)
\end{equation}
where $h_n$ is the timestep, $\jac_n = f'(u_n)$ is the Jacobian matrix, and $u_{n}$ is the numerical approximation to $u(t_{n})$.
The matrix function $\varphi_1(z) = (\exp(z) - 1)/z$ is just one member of the sequence of functions
\begin{equation} \label{eq:phi}
    \varphi_{0}(z) = \exp(z),
    \qquad
    \varphi_{k+1}(z) = \frac{\varphi_k(z) - \varphi_k(0)}{z},
\end{equation}
which are ubiquitous in the exponential integrator literature.

It is well-known that the exponential Rosenbrock--Euler method is second order accurate both in the classical sense \cite{pope1963exponential} and for stiff, semilinear problems \cite{hochbruck2009exponential}.
These results are based on the assumption that $\varphi_1$ is computed exactly; however, this is rarely the case in practice.
Typically, it is computed to a specified tolerance and contains errors from floating point arithmetic.
As demonstrated in experiments later in this section, performing the linear algebra associated with $\varphi$-functions on low precision hardware can severely limit the accuracy of an exponential integrator.

Using \eqref{eq:phi}, we can equivalently express the exponential Rosenbrock--Euler scheme \eqref{eq:exp_Euler} as
\begin{equation} \label{eq:reformulation}
    \begin{split}
        u_{n+1} &= u_n + h_n f(h_n) + h_n (\varphi_1(h_n \jac_n) - I) f(u_n) \\
        &= u_n + h_n f(u_n) + h_n^2 \varphi_2(h_n \jac_n) \jac_n f(u_n).
    \end{split}
\end{equation}
The benefit of using formulation \eqref{eq:reformulation} is the $\varphi$-function is scaled by $h_n^2$ as opposed to $h_n$ in \eqref{eq:exp_Euler}.
Consequently, we may expect improved resilience to $\varphi$-function errors as $h_n \to 0$.
However, this asymptotic analysis breaks down when $\jac_n$ is disproportionally large and $h_n$ is not sufficiently small.
In this stiff regime, the term $h^2 \varphi_2(h \jac_n) \jac_n f(u_n)$ is susceptible to overflows as well as cancellation errors with $h_n f(u_n)$.

Therefore, we propose the following reformulated exponential Rosenbrock--Euler scheme which uses a parameter, $\gamma_n$, to vary between \eqref{eq:exp_Euler} and \eqref{eq:reformulation},
\begin{equation} \label{eq:mixed_precision_exp_Euler}
    \ump_{n+1} = \ump_n + h_n \gamma_n f(\ump_n) + \lp(h_n \psi(h_n \jacmp_n, \gamma_n) f(\ump_n)).
\end{equation}
We use $\ump_n$ and $\jacmp_n = f'(\ump_n)$ to distinguish  from $u_n$ which uses exact $\varphi$-functions.
The function $\lp(x)$ represents the evaluation of $x$ to a tolerance $\epsilon$ and is assumed to satisfy the error model $\lp(x) = (I + \delta) x$ with $\norm{\delta}_2 \leq \epsilon$.
Finally, we introduce
\begin{subequations}
    \begin{align}
        \label{eq:psi:1}
        \psi(z, \gamma) &= \varphi_1(z) - \gamma \\
        \label{eq:psi:2}
        &= (1 - \gamma) \varphi_1(z) + \gamma \varphi_2(z) z.
    \end{align}
\end{subequations}
While the form \eqref{eq:psi:1} is useful for analysis, \eqref{eq:psi:2} is preferable for implementation as it is less susceptible to subtractive cancellation.

\subsection{Error Analysis}

In order to inform the selection of the yet unspecified parameter $\gamma_n$ in \eqref{eq:mixed_precision_exp_Euler}, we first study the effect of $\gamma_n$ on the numerical error.
The local truncation error committed after one step is
\begin{equation*}
    e_1 = \ump_1 - u(t_1).
\end{equation*}
This satisfies
\begin{equation}
    \begin{split} \label{eq:error}
        \norm{e_1}_2 &= \norm{\ump_1 - u(t_1)}_2 \\
        & \leq \norm{u_1 - u(t_1)}_2 + \norm{\ump_1 - u_1}_2 \\
        & \leq C h_0^3 + \norm{\delta h_0 \psi(h_0 \jacmp_0, \gamma_0) f(\ump_0)}_2 \\
        & \leq C h_0^3 + h_0 \epsilon \norm{\varphi_1(h_0 \jacmp_0) f(\ump_0) - \gamma_0 f(\ump_0)}_2,
    \end{split}
\end{equation}
where we have used the triangle inequality and the second order convergence property of the exponential Rosenbrock--Euler method.

This suggests solving the optimization problem
\begin{equation} \label{eq:gamma_opt}
    \begin{split}
        \gamma_n &= \argmin_{\gamma} \norm{\varphi_1(h_n \jacmp_n) f(\ump_n) - \gamma f(\ump_n)}^2_2 \\
        &= \frac{f(\ump_n)^T \varphi_1(h_n \jacmp_n) f(\ump_n)}{\norm{f(\ump_n)}^2_2}
    \end{split}
\end{equation}
to select $\gamma_n$ at each step to minimize the effect of the low precision arithmetic.
As $\varphi_1$ is already required to compute $\psi$ in \eqref{eq:mixed_precision_exp_Euler}, the additional cost of computing $\gamma_n$ is negligible for many algorithms used to compute linear combinations of $\varphi$-functions.
Alternatively, one can use the bound
\begin{equation} \label{eq:gamma_est}
    %\frac{f(\ump_n)^* \varphi_1(h_n \jacmp_n) f(\ump_n)}{\norm{f(\ump_n)}^2_2}
    \gamma_n
    \leq \mu_2(\varphi_1(h_n \jacmp_n))
    \leq \varphi_1(h_n \mu_2(\jacmp_n))
\end{equation}
to choose $\gamma_n$.
If $\mu_2(\jacmp_n)$, the logarithmic 2-norm \cite{soderlind2006logarithmic} of the Jacobian, can be readily estimated, \eqref{eq:gamma_est} only requires inexpensive scalar arithmetic.

In the stiff regime where $\jacmp_n \to -\infty$, 
%\note[id=DJG]{Should this be the norm of $\jacmp_n$?}
%\note[id=SBR]{A norm is insufficient because $f$ may be orthogonal to eigenvectors of diverging eigenvalues. I couldn't think of simple matrix condition, so I used this scalar one. Open to ideas} 
the reformulated method \eqref{eq:mixed_precision_exp_Euler} approaches \eqref{eq:exp_Euler} because $\gamma_n \to 0$.
Conversely, in the asymptotic regime where $h_n \to 0$, \eqref{eq:mixed_precision_exp_Euler} approaches \eqref{eq:reformulation}.
A Taylor expansion of \eqref{eq:gamma_opt} reveals $\gamma_n = 1 - \order{h}$.
Thus, the local truncation error is $\norm{e_1}_2 = \order{h_0^3 + \epsilon h_0^2}$ as opposed to $\order{h_0^3 + \epsilon h_0}$ for the traditional formulation \eqref{eq:exp_Euler}.

\subsection{Convergence Experiment}

In order to verify the improved accuracy of \eqref{eq:mixed_precision_exp_Euler}, we compare its convergence to \eqref{eq:exp_Euler} on an advection–diffusion–reaction PDE from \cite[Section 5.1]{caliari2009implementation},
\begin{equation}\label{eq:adr_pde}
    \begin{split}
        \pdv{u}{t} &= \varepsilon \left( \pdv[2]{u}{x} + \pdv[2]{u}{y} \right) - \alpha \left( \pdv{u}{x} + \pdv{u}{y} \right) \\
        & \quad + \rho u \left( u - \frac{1}{2} \right) (1 - u), \\
        u(0, x, y) &= 0.3 + 256 (x(1-x)y(1-y))^2.
    \end{split}
\end{equation}
The timespan is $[0, 0.3]$ and the spatial domain, $x,y \in [0,1]$, is discretized by second order finite differences with $\Delta x = \Delta y = 0.05$.
The remaining parameters are $\varepsilon = 0.05$, $\alpha = -1$, and $\rho = 1$.
Error in the numerical solution is measured as $\norm{u_n - u_\text{ref}}_2$, where $u_\text{ref}$ is a reference solution computed with an absolute and relative tolerance of $10^{-13}$.

\begin{figure}[ht!]
    \centering
    \begin{tikzpicture}
        \begin{loglogaxis}[xlabel={Time Steps}, ylabel={Absolute $\ell^2$ Error}, legend entries={Exp Ros--Euler, Reformulated Exp Ros--Euler}, legend pos=south west, max space between ticks=20, width=0.95\linewidth, height=0.55\linewidth]
            \addplot table [x index=0, y index=1] {rda_21x21_convergence_f4.csv};
            \addplot table [x index=0, y index=2] {rda_21x21_convergence_f4.csv};
            \draw (axis cs:1e2,3e-6) -- node[below]{$2$} (axis cs:1e3,3e-8);
        \end{loglogaxis}
    \end{tikzpicture}
    
    \caption{The reformulated exponential Rosenbrock--Euler method \eqref{eq:mixed_precision_exp_Euler} maintains second order convergence despite using single precision for $\varphi$-functions, while the traditional form \eqref{eq:exp_Euler} stagnates at the accuracy of the $\varphi$-function.}
    \label{fig:rda_convergence_f4}
\end{figure}
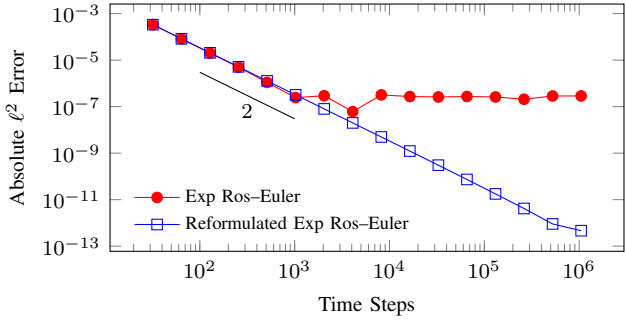

Our first experiment uses single precision for the terms in the $\lp$ function of \eqref{eq:mixed_precision_exp_Euler} including the Jacobian evaluation and $\varphi$-functions.
The remaining operations, including evaluating $f$, are performed in double precision.
Figure \ref{fig:rda_convergence_f4} shows that the accuracy of the traditional exponential Rosenbrock--Euler form is limited by the accuracy of the $\varphi$-functions as it cannot achieve an error below $10^{-7}$.
The reformulated version \eqref{eq:mixed_precision_exp_Euler}, however, is able to achieve errors six orders of magnitude smaller without suffering order reduction.

\begin{figure}[ht!]
    \centering
    \begin{tikzpicture}
        \begin{loglogaxis}[xlabel={Time Steps}, ylabel={Absolute $\ell^2$ Error}, legend entries={Exp Ros--Euler, Reformulated Exp Ros--Euler}, legend style={at={(0.55,0.98)},anchor=north}, max space between ticks=20, width=0.95\linewidth, height=0.55\linewidth]
            \addplot table [x index=0, y index=1] {rda_21x21_convergence_f2.csv};
            \addplot table [x index=0, y index=2] {rda_21x21_convergence_f2.csv};
            \draw (axis cs:2.4e0,2e-2) -- node[left]{$2$} (axis cs:8e0,1.5e-3);
            \draw (axis cs:1.5e4,1.6e-4) -- node[below]{$1$} (axis cs:6e4,4e-5);
        \end{loglogaxis}
    \end{tikzpicture}
    
    \caption{With half precision $\varphi$-functions, the reformulated exponential Rosenbrock--Euler method \eqref{eq:mixed_precision_exp_Euler} achieves a minimum error approximately ten times smaller than that of the traditional form \eqref{eq:exp_Euler}.}
    \label{fig:rda_convergence_f2}
\end{figure}
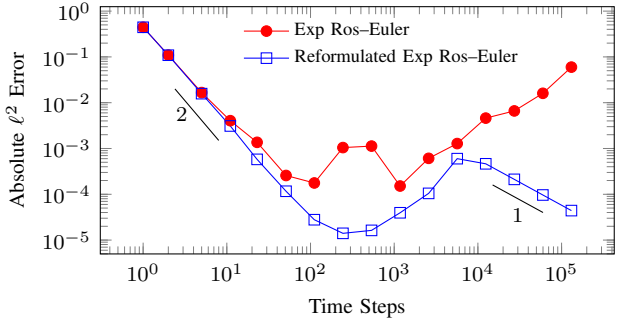

When using half precision instead of single precision for the Jacobian evaluation and $\varphi$-functions (Figure \ref{fig:rda_convergence_f2}), the improvement in accuracy with the reformulated method \eqref{eq:mixed_precision_exp_Euler} is more modest: an order of magnitude. After a momentary degradation in accuracy between $200$ and $6000$ time steps, the $\epsilon h_0^2$ term of the local error becomes dominant and we see further asymptotic improvements.
%\note[id=DJG]{If I remeber correctly, similar behavior was observed in Grant's ``Perturbed Runge–Kutta Methods for Mixed Precision Applications'' paper.}

\section{Computing \texorpdfstring{$\varphi$}{phi}-function products with a mixed precision Krylov method}
\label{s:KIIOPS}

Historically, exponential methods were bound by the cost and difficulty of computing matrix exponentials in the $\varphi$-functions. However, over the last several decades a rich literature has developed around the use of Krylov approximations for the action of $\varphi$-functions on a vector \cite{sidje1998expokit,niesen2012algorithm,almohy2011computing,gaudreault2018kiops}. 
These Krylov-based approaches have made exponential integrators more practical to use and are based on the approximation,
\begin{equation}\label{eq:krylov_approximation}
\begin{split}
    & \quad \exp(\tau A) v = \beta V_m \exp(\tau H_m) e_1 , \\
    & m \ll n,\quad \beta = \| v \|_2,\quad e_1 = (1, 0, \dots, 0)^T,
\end{split}
\end{equation}
where $V_m \in \mathbb{R}^{n \times m}$ is the orthonormal basis of the Krylov subspace $\mathcal{K}_m(A, v)$ and $H_m \in \mathbb{R}^{m \times m}$ is the Hessenberg matrix generated by the Arnoldi process. 
In early approaches for computing \eqref{eq:krylov_approximation}, such as the methods in Expokit \cite{sidje1998expokit}, the computational cost is dominated by the full orthogonlization method (FOM) \cite{voSidje2017}. 
State-of-the-art implementations, such as \KIOPS{} \cite{gaudreault2018kiops}, utilize an incomplete orthogonalization process (IOP).
Using IOP is not only faster, but it shifts the majority of the computing effort to matrix-vector products \cite{voSidje2017} which map well to low precision computing units on modern hardware
% Using IOP is not only faster, but opens the door to leveraging low precision to accelerate the computations because it is less susceptible roundoff accumulation than FOM [CITE, paper that discusses roundoff in FOM] and shifts the majority of the computing effort to matrix-vector products \cite{voSidje2017} which map well to low precision computing units on modern hardware.

\subsection{Introducing low precision into IOP Arnoldi}

It is therefore natural to consider introducing low precision computations into IOP Arnoldi via the matrix-vector products since they map well to many low precision hardware units. 
A naive approach is simply to perform all of the matrix-vector products in low precision. 
We demonstrate the problems with this approach by modifying the IOP Arnoldi procedure (\ref{alg:IOP_Arnoldi}) in \KIOPS{} so that the matrix-vector products are computed in low precision and the result is stored in double-precision (we use the \texttt{chop} method from \cite{higham2019simulating} to simulate this) in two different experiments with IEEE single, NVIDIA TensorFloat-32 (TF32), IEEE half, and bfloat16 floating-point formats.
The test setups we use are essentially the same as the ones utilized by Al-Mohy and Higham \cite[Experiment 5 and 7]{almohy2011computing} which are based on experiments conducted by Niesen and Wright \cite[Experiment 1]{niesen2012algorithm} and Sidje \cite[Section 6.2]{sidje1998expokit}.
% We note that during the course of investigation Caliari et al. introduced the {\small BAMPHI} method (based on Newton interpolation at the extended Ritz values of $A$) \cite{caliari2022bamphi} that seems to outperform \KIOPS{}. {\small BAMPHI} also employs an IOP Arnoldi procedure, thus we suspect the mixed precision methods described in this section may be adaptable to this approach as well. However, we do not extend our experiments in this paper.

\subsubsection*{Experiment 1}

Compute $u = \exp(tA)b_0$ using \KIOPS{} with Algorithm \ref{alg:IOP_Arnoldi} for three matrices.
The first two matrices are from the Harwell-Boeing collection \cite{duff1989sparse} and are available in the SuiteSparse Sparse Matrix Collection \cite{davis2011university}. 
With the \texttt{orani678} sparse matrix (order $n = 2529$ with $nnz = 90158$ nonzero elements) we use $t = 10$, $b_0 = [1, \dots, 1]^T$, and set the \KIOPS{} tolerance to $tol = \sqrt{\epsilon_d}$ where $\epsilon_d$ is machine precision for the double precision format. 
For the \texttt{bcspwr10} sparse matrix (order $n = 5300$ with $nnz = 21842$) we use $t = 10$, $b_0 = [1, 0, \dots, 0, 1]^T$, and $tol = 10^{-5}$.
The third test uses a Poisson matrix of order $n = 9801$, $t = 1$, and $tol = 10^{-12}$. The matrix and the $b$ vector are generated with the MATLAB code
\begin{verbatim}
A = 2500 * gallery('poisson', 99);
g = (-0.98 : 0.02 : 0.98)';
[R1, R2] = meshgrid(g, g); 
r1 = R1(:); r2 = R2(:);
b = (1 - r1.^2) .* (1 - r2.^2) .* exp(r1);
\end{verbatim}

% The second matrix is the \texttt{bcspwr10} sparse matrix of order $n = 5300$ with $nnz = 21842$ nonzero elements. 
% With the \texttt{bcspwr10} matrix we use $t = 2$ and $b_0 = [1, 0, \dots, 0, 1]^T$, and $tol = 10^{-5}$. \note[id=DJG]{Why the difference in tolerance?}\note[id=SBR]{Let's include the value of $\sqrt{\epsilon}$} 

\subsubsection*{Experiment 2}

Compute $u = \varphi_0(tA)b_0 + t \varphi_1(tA)b_1 + \cdots + t^4 \varphi_4(tA)b_4$ via the modified \KIOPS{} method with the \texttt{orani678}, \texttt{bcspwr10}, and Poisson matrices and $b_i = [1, \dots, 1]^T$ using the same values for $t$ and $tol$ as in Experiment 1.

\begin{algorithm}
  \caption{Naive low precision IOP Arnoldi with the low precision computation (line \ref{line:IOP_arnoldi_mv1}) boxed.}\label{alg:IOP_Arnoldi}
  \begin{algorithmic}[1]
  \STATE \textbf{Input:} $A \in \mathbb{R}^{N \times N}$, $B \in \mathbb{R}^{N \times p}$, $V \in \mathbb{R}^{(N + p) \times (m_{\max} + 1)}$, $j$, $m$
  \WHILE { $j < m$ }
     \STATE $j = j + 1$
     \STATE $V(1:N, j+1) =$\par\qquad$ \boxed{A \cdot V(1:N, j)  + B \cdot V(N+1:N+p, j)}$\label{line:IOP_arnoldi_mv1}
     \STATE $V(N+1:N+p-1, j+1) = V(N+2:N+p, j)$
     \STATE $V(N+p, j+1) = 0$
     \FOR { $i =  \max(1,j-1)$ \TO $j$ }
        \STATE   $H(i,j) = V(:,i)^{T} \cdot V(:, j+1)$
        \STATE   $V(:, j + 1) = V(:, j + 1) - H(i, j) \cdot V(:, i)$
     \ENDFOR
  
     \STATE $s = \| V(:, j + 1) \|_2$
  
     \IF { $s \approx 0$ }
        \STATE happy\_breakdown = true
        \STATE \textbf{break}
     \ENDIF
  
     \STATE $H(i+1,j) = s$
     \STATE $V(:,j+1) = \frac{1}{s} V(:, j + 1) $
  \ENDWHILE
  \RETURN $V$, $H$, $j$
  \end{algorithmic}
\end{algorithm}

\begin{table*}[!htb]
  \centering
  \begin{tabular}{|c|cc|c|c|c|c|c|}
      \hline
      &  &  & double & single & TF32 & half & bfloat16 \\
      \hline
      & matrix & $tol$ & error & error & error & error & error \\
      \hline
      \rule{0pt}{3ex} 
      \multirow{3}{*}{Experiment 1} & 
        orani678  & $1.49\mathrm{e}{-8}$  & $6.88\mathrm{e}{-12}$ & $2.41\mathrm{e}{-5}$ & $1.32\mathrm{e}{-1}$  & $1.71\mathrm{e}{+0}$ & $3.91\mathrm{e}{+1}$ \\
      & bcspwr10  & $1.00\mathrm{e}{-5}$  & $4.84\mathrm{e}{-10}$ & $3.10\mathrm{e}{-2}$ & $4.03\mathrm{e}{+2}$  & $7.67\mathrm{e}{+3}$ & $3.39\mathrm{e}{+3}$ \\
      & Poisson   & $1.00\mathrm{e}{-12}$ & $8.38\mathrm{e}{-14}$ & $1.15\mathrm{e}{-9}$ & $1.76\mathrm{e}{-5}$  & $1.14\mathrm{e}{-5}$ & $7.70\mathrm{e}{-3}$ \\
      \hline
      \rule{0pt}{3ex} 
      \multirow{3}{*}{Experiment 2} & 
        orani678  & $1.49\mathrm{e}{-8}$  & $1.69\mathrm{e}{-13}$ & $2.91\mathrm{e}{-5}$ & $7.63\mathrm{e}{-1}$ & $1.09\mathrm{e}{+0}$ & $2.87\mathrm{e}{+0}$ \\
      & bcspwr10  & $1.00\mathrm{e}{-5}$  & $1.73\mathrm{e}{-5}$ & $2.37\mathrm{e}{-6}$ & $2.06\mathrm{e}{-2}$ & $2.06\mathrm{e}{-2}$ & $1.27\mathrm{e}{-1}$ \\
      & Poisson   & $1.00\mathrm{e}{-12}$ & $2.06\mathrm{e}{-14}$ & $1.26\mathrm{e}{-9}$ & $4.60\mathrm{e}{-3}$ & $6.40\mathrm{e}{-3}$ & $1.85\mathrm{e}{-1}$ \\
      \hline
  \end{tabular}
  \vspace{1mm}
  \caption{Results from experiments with the naive low precision IOP Arnoldi (Algorithm \ref{alg:IOP_Arnoldi}) in \KIOPS{} show that the method is unreliable. In all cases, the relative error is large when using low precision with respect to both the tolerance and the relative error achieved with double precision.}
  \label{tab:naive_mixed_iop_arnoldi_1}
\end{table*}

Letting $u_p$ be the solution generated using \KIOPS{} with Algorithm \ref{alg:IOP_Arnoldi} and precision $p$, we define the error as  $err(u_p) = \| u_p - u_\text{ref} \|_{\infty}/\|u_\text{ref}\|_{\infty}$.
The reference solution $u_\text{ref}$ is generated with the standard \KIOPS{} method in double precision with a tolerance of $tol = \epsilon_d$. Unless otherwise stated, results use the default \KIOPS{} parameters.
When using the naively modified IOP Arnoldi in \KIOPS{}, we see that the error is far greater than the desired tolerance (Table \ref{tab:naive_mixed_iop_arnoldi_1}).

In an attempt to recover the lost accuracy from low precision matrix-vector products, we now reconsider replacing the exact (in finite arithmetic) matrix-vector products with the inexact matrix-vector product
\begin{equation}\label{eq:inexact_product}
  \tilde{A} v = (A + E)v,
\end{equation}
where $E$ is some perturbation matrix.
Substituting \eqref{eq:inexact_product} into \eqref{eq:krylov_approximation} and allowing the $E$ to change with the Arndoli iterate yields the inexact Arnoldi approximation
\begin{equation}\label{eq:inexact_Arnoldi}
\begin{split}
  (A + \fancyE_m)V_m &= V_m H_m + h_{m+1,m} v_{m+1} e^T_m, \\
  \fancyE_m &= \sum_{j=1}^m E_j v_j v^T_j.
\end{split}
\end{equation}
The theoretical underpinnings for this approach are developed in \cite{simoncini2003Inexact}. Furthermore, \cite{bouras2005Inexact} provides bounds on the growth of $\|E_j\|_2$ as the iterations progress in various Krylov subspace methods including FOM Arnoldi.
Dinh and Sidje  extended the work to computing the matrix-exponential in \cite{dinh2017Analysis}.
However, the combination of IOP, inexact products, and matrix-exponential computations has, as far as we are aware, not been previously examined in the literature.

%A full mathematical analysis of this combination is outside the scope of this work 
We numerically investigate the effectiveness of this intuitive approach by progressively introducing lower-precision matrix-vector products (i.e., allowing $\| E_j \|_2$ to grow) into the IOP Arnoldi algorithm within \KIOPS{} as the Arnoldi iteration proceeds.
We define two new parameters $\mchop{1}$ and $\mchop{2}$ that determine the Arnoldi iterates at which we switch from full double-precision matrix-vector products to single-precision and then from single to either TF32, half, or bfloat16 (Algorithm \ref{alg:IOP_Arnoldi_Smart}).
% That is, when $j$ in Algorithm \ref{alg:IOP_Arnoldi_Smart} is greater than $\mchop{1}$ we switch to using single precision for the matrix-vector products, and when $j$ is greater than $\mchop{2}$ we switch to using half, TF32, or bfloat16. 

We repeat Experiments 1 and 2 while first varying $\mchop{1}$ until the target error is below $\max(err(u_d), tol)$.
This metric is employed because $err(u_d) > tol$ in Experiment 2 with the \texttt{bcspwr10} matrix, so there is no hope of doing better than $err(u_d)$ in this case.
Then, with $\mchop{1}$ fixed to the value we just found, we vary $\mchop{2}$ until the tolerance is met.
Utilizing this procedure with mixed precision IOP Arnoldi enables \KIOPS{} to achieve a much lower error while leveraging a precision lower than double for 40\% or more of the Arnoldi iterates (Table \ref{tab:smarter_mixed_iop_arnoldi}). 
Furthermore, we are able to leverage lower than single-precision for 25\% -- 60\% of iterates.
% Given the significantly increased FLOP rates (((refer to Table with hardware))) for these low precision matrix-vector products, Algorithm \ref{alg:IOP_Arnoldi_Smart} clearly has the potential to provide a significant speedup to exponential time-integrators.
% \textcolor{blue}{It would be nice to have a result showing the speedup. If we cannot produce one in time, we should porbably scrap the previous sentence or at least say its future work.}

\begin{table*}[!htb]
    \centering
    \tabcolsep=0.15cm
    \begin{tabular}{|c|cc|ccc|ccc|ccc|}
        \hline
        &        &     & \multicolumn{3}{|c|}{single, TF32} & \multicolumn{3}{|c|}{single, half} & \multicolumn{3}{|c|}{single, bfloat16} \\
        \hline
        & matrix & target error & $\mchop{1,2}$ & $m$ & error & $\mchop{1,2}$ & $m$ & error & $\mchop{1,2}$ & $m$ & error  \\
        \hline
        \rule{0pt}{3ex} 
        \multirow{3}{*}{Experiment 1} & 
          orani678  & $1.49\mathrm{e}{-8}$  & 30, 39 & 51  & $4.01\mathrm{e}{-9}$  & 30, 39 & 51  & $4.01\mathrm{e}{-9}$  & 30, 40 & 51  & $5.41\mathrm{e}{-9}$  \\
        & bcspwr10  & $1.00\mathrm{e}{-5}$  & 36, 54 & 86  & $6.82\mathrm{e}{-6}$  & 36, 54 & 86  & $6.82\mathrm{e}{-6}$  & 36, 54 & 87  & $6.64\mathrm{e}{-6}$  \\
        & Poisson   & $1.00\mathrm{e}{-12}$ & 15, 70 & 128 & $6.49\mathrm{e}{-13}$ & 15, 70 & 128 & $6.68\mathrm{e}{-13}$ & 15, 90 & 128 & $9.64\mathrm{e}{-13}$ \\
        \hline
        \rule{0pt}{3ex} 
        \multirow{3}{*}{Experiment 2} & 
          orani678  & $1.49\mathrm{e}{-8}$  & 28, 37 & 51  & $2.85\mathrm{e}{-9}$  & 28, 37 & 51  & $2.85\mathrm{e}{-9}$  & 28, 39 & 51  & $2.61\mathrm{e}{-9}$  \\
        & bcspwr10  & $1.73\mathrm{e}{-5}$  & 51, 63 & 106 & $1.73\mathrm{e}{-5}$  & 51, 63 & 106 & $1.73\mathrm{e}{-5}$  & 51, 66 & 106 & $1.73\mathrm{e}{-5}$  \\
        & Poisson   & $1.00\mathrm{e}{-12}$ & 28, 58 & 128 & $8.89\mathrm{e}{-13}$ & 28, 58 & 128 & $9.26\mathrm{e}{-13}$ & 28, 60 & 128 & $3.40\mathrm{e}{-13}$ \\
        \hline
    \end{tabular}
    \vspace{1mm}
    \caption{Progressively introducing lower-precision matrix-vector products into the IOP Arnoldi procedure (Algorithm \ref{alg:IOP_Arnoldi_Smart}) within \KIOPS{} enables the target error, $\max(err(u_d), tol)$, to be met. $\mchop{1}$ and $\mchop{2}$ are the Krylov iterations for switching from single-precision to TF32, half, or bfloat16. $m$ is the number of vectors in $\mathcal{K}(A,v)$ for the last iteration of \KIOPS{} and is generally a good estimate for the basis size required (128 is the default maximum).}
    \label{tab:smarter_mixed_iop_arnoldi}
\end{table*}

\begin{algorithm}
  \caption{Mixed precision IOP Arnoldi with the first low precision computations (line \ref{line:mv1}) in the dashed box and lowest precision the in solid box (line \ref{line:mv2}).}\label{alg:IOP_Arnoldi_Smart}
  \begin{algorithmic}[1]
  \STATE \textbf{Input:} $A \in \mathbb{R}^{N \times N}$, $B \in \mathbb{R}^{N \times p}$, $V \in \mathbb{R}^{(N + p) \times (m_{\max} + 1)}$, $j$, $m$
  \WHILE { $j < m$ }
     \STATE $j = j + 1$
     \IF { $j + 1 > m_{\text{chop}_2}$ }
        \STATE $V(1:N, j+1) =$\par\qquad$ \boxed{A \cdot V(1:N, j)  + B \cdot V(N+1:N+p, j)}$\label{line:mv2}
      \ELSIF { $j + 1 > m_{\text{chop}_1}$ }
        \STATE $V(1:N, j+1) =$\par\qquad \dbox{$A \cdot V(1:N, j)  + B \cdot V(N+1:N+p, j)$}\label{line:mv1}
     \ELSE 
        \STATE $V(1:N, j+1) =$\par\qquad$ A \cdot V(1:N, j)  + B \cdot V(N+1:N+p, j)$
     \ENDIF
     \STATE $V(N+1:N+p-1, j+1) = V(N+2:N+p, j)$
     \STATE $V(N+p, j+1) = 0$ 
     \FOR { $i =  \max(1,j-1)$ \TO $j$ }
        \STATE   $H(i,j) = V(:,i)^{T} \cdot V(:, j+1)$
        \STATE   $V(:, j + 1) = V(:, j + 1) - H(i, j) \cdot V(:, i)$
     \ENDFOR
  
     \STATE $s = \| V(:, j + 1) \|_2$
  
     \IF { $s \approx 0$ }
        \STATE happy\_breakdown = true
        \STATE \textbf{break}
     \ENDIF
  
     \STATE $H(i+1,j) = s$
     \STATE $V(:,j+1) = \frac{1}{s} V(:, j + 1) $
  \ENDWHILE
  \RETURN $V$, $H$, $j$
  \end{algorithmic}
\end{algorithm}

\section{Integrated Numerical Experiments}
\label{s:tests}

To evaluate the performance of the two approaches for leveraging low precision computation in exponential integrators we test three methods: standard exponential Rosenbrock--Euler, the reformulated exponential Rosenbrock--Euler scheme \eqref{eq:mixed_precision_exp_Euler}, and the stiffly-accurate fourth-order exprk4s6 \cite{luan2021efficient}.
All three methods are tested with standard \KIOPS{} and \KIOPS{} with the mixed precision IOP Arnoldi (Algorithm \ref{alg:IOP_Arnoldi_Smart}) for evaluating $\varphi$-function vector products. As before we use \texttt{chop} for simulating low-precision computations. 
With mixed precision IOP Arnoldi we use $m_{\text{chop}_1}$ and $m_{\text{chop}_2}$ to set the iteration for switching to single or half precision matrix-vector products, respectively.
The process used to choose these values is similar to the process used in the experiments in Section \ref{s:KIIOPS}. 
We find a value for $m_{\text{chop}_{1}}$ that produces the desired error, fix its value, then we find $m_{\text{chop}_{2}}$ that similarly allows the desired error to be met. 
For the exprk4s6 method, this means we have to choose the values for each of the four calls per time step that it makes to \KIOPS{}.
The six possible combinations of schemes are used to solve the advection-diffusion-reaction problem \eqref{eq:adr_pde} with the same parameters but a finer spatial discretization, $\Delta x = \Delta y = 0.0025$, leading to a stiffer problem.
% \note[id=DR]{You say that $\Delta x = \Delta y = 0.0025$ leads to a “stiffer problem” — it would help if you would *quantify* how much stiffer it really is (e.g., $\lambda_{max}(J) / \lambda_{min}(J)$).}

\begin{figure*}[!htb]
    \centering
    \begin{tikzpicture}
        \begin{groupplot}[group style={group size= 3 by 1, ylabels at=edge left, horizontal sep=1.1cm}, xlabel={Time Steps}, ylabel={Relative $l_{\infty}$ Error}, xmode=log, ymode=log, width=0.35\textwidth, legend pos=south west]
            \nextgroupplot[legend entries={ERE KIOPS-double, ERE KIOPS-half, ERE KIOPS-mixed}]
            \addplot table [x=steps, y=ere_dbl_error] {rda_400x400_convergence_exp_rbeuler.csv};
            \addplot table [x=steps, y=ere_half_error] {rda_400x400_convergence_exp_rbeuler.csv};
            \addplot table [x=steps, y=ere_mixed_error] {rda_400x400_convergence_exp_rbeuler.csv};
            \nextgroupplot[legend entries={ERE KIOPS-double, RERE KIOPS-half, RERE KIOPS-mixed}]
            \addplot table [x=steps, y=ere_dbl_error] {rda_400x400_convergence_exp_rbeuler.csv};
            \addplot table [x=steps, y=rere_half_error] {rda_400x400_convergence_exp_rbeuler.csv};
            \addplot table [x=steps, y=rere_mixed_error] {rda_400x400_convergence_exp_rbeuler.csv};
            \nextgroupplot[legend entries={exprk4s6 KIOPS-double, exprk4s6 KIOPS-mixed}]
            \addplot table [x=steps, y=dbl_error] {rda_400x400_convergence_exprk4s6.csv};
            \addplot table [x=steps, y=mixed_error] {rda_400x400_convergence_exprk4s6.csv};
        \end{groupplot}
    \end{tikzpicture}
    
    \caption{The reformulated exponential Rosenbrock--Euler (RERE) method \eqref{eq:mixed_precision_exp_Euler} achieves a much lower error than the standard exponential Rosenbrock--Euler (ERE) when using \KIOPS{} with only half precision matrix-vector products via Algorithm \ref{alg:IOP_Arnoldi}. Using mixed precision matrix-vector products in \KIOPS{} via Algorithm \ref{alg:IOP_Arnoldi_Smart} significantly improves the error for all methods. exprk4s6 is unable to converge at all with Algorithm \ref{alg:IOP_Arnoldi} and only half precision matrix-vector products while using Algorithm \ref{alg:IOP_Arnoldi_Smart} enables it to run and obtain a reasonably accurate solution.}
    \label{fig:order_plots}
\end{figure*}
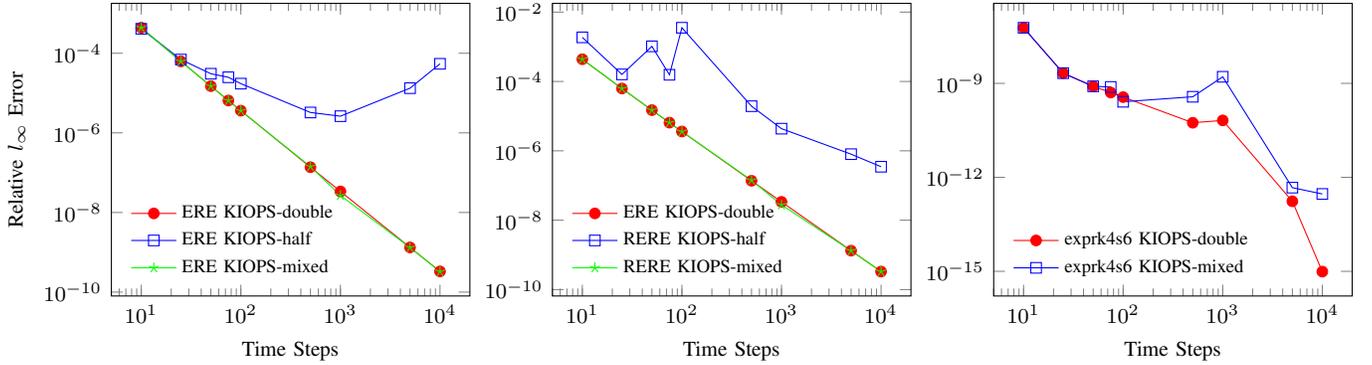

\subsection{Overall accuracy}

Figure \ref{fig:order_plots} shows the error $\| u_n - u_\text{ref} \|_{\infty} / \| u_\text{ref} \|_{\infty}$ versus the number of time steps to demonstrate the convergence of the different schemes. 
The reference solution is generated with the exprk4s6 method with $10^5$ time steps.
Once again we see that the reformulated exponential Rosenbrock--Euler method \eqref{eq:mixed_precision_exp_Euler} consistently achieves lower error and maintains second order convergence longer than the standard Rosenbrock--Euler method.
The use of \KIOPS{} with mixed precision IOP Arnoldi greatly improves the accuracy for both methods, with the error nearly identical to what is achieved when using \KIOPS{} with double precision.
In the case of the higher-order exprk4s6, mixed precision IOP Arnoldi enables using low precision as running with only half precision does not converge.
% \added{With the standard Rosenbrock--Euler method we use $m_{\text{chop}_{1,2}} = 2, 3$, but with the mixed precision exponential Rosenbrock--Euler method we use $m_{\text{chop}_{1,2}} = 1, 3$ as its additional robustness to low precision enables utilizing single precision sooner.} \note[id=DJG]{Could we combine these plots into one (or two, standard and mixed precision KIOPS) plots}

\subsection{Idealized computational efficiency}

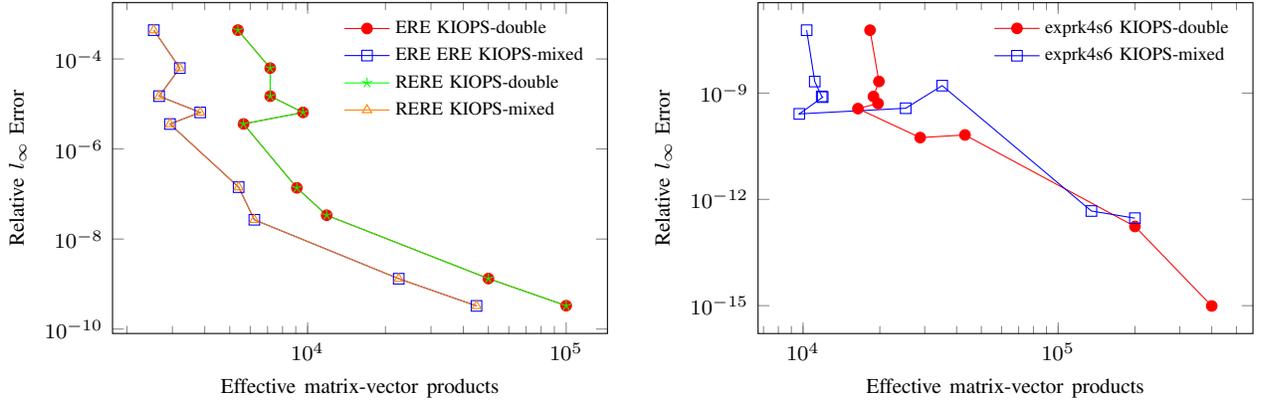
\begin{figure*}[!htb]
    \centering
    \begin{tikzpicture}
        \begin{groupplot}[group style={group size= 2 by 1, horizontal sep=2cm}, xlabel={Effective matrix-vector products}, ylabel={Relative $l_{\infty}$ Error}, xmode=log, ymode=log, width=0.45\textwidth, height=0.33\textwidth]
            \nextgroupplot[legend entries={ERE KIOPS-double, ERE ERE KIOPS-mixed, RERE KIOPS-double, RERE KIOPS-mixed}]
            \addplot table [y=ere_dbl_err, x=ere_dbl_mv] {rda_400x400_work_exp_rbeuler.csv};
            \addplot table [y=ere_mixed_err, x=ere_mixed_mv] {rda_400x400_work_exp_rbeuler.csv};
            \addplot table [y=rere_dbl_err, x=rere_dbl_mv] {rda_400x400_work_exp_rbeuler.csv};
            \addplot table [y=rere_mixed_err, x=rere_mixed_mv] {rda_400x400_work_exp_rbeuler.csv};
            \nextgroupplot[legend entries={exprk4s6 KIOPS-double, exprk4s6 KIOPS-mixed}]
            \addplot table [y=exprk_dbl_err, x=exprk_dbl_mv] {rda_400x400_work_exp_rbeuler.csv};
            \addplot table [y=exprk_mixed_err, x=exprk_mixed_mv] {rda_400x400_work_exp_rbeuler.csv};
        \end{groupplot}
    \end{tikzpicture}
    
    \caption{In the case of the reformulated exponential Rosenbrock--Euler (RERE) and standard exponential Rosenbrock--Euler (ERE) methods, using the mixed precision IOP Arnoldi (Algorithm \ref{alg:IOP_Arnoldi_Smart}) in \KIOPS{} reduced the number of effective matrix-vector products compared to the double precision version. With exprk4s6, the mixed precision \KIOPS{} was more efficient until the error reached $~10^{-10}$.}
    \label{fig:work_precision}
\end{figure*}

Figure \ref{fig:work_precision} provides an estimate of the computational efficiency and shows the error versus the number of ``effective'' matrix-vector products, $mv_{\text{effective}}$, where
\begin{equation}
    mv_{\text{effective}} = mv_{\text{double}} + \frac{mv_{\text{single}}}{a} + \frac{mv_{\text{half}}}{b}.
\end{equation}
We use matrix-vector products as a proxy for the wall-clock time since they are typically the critical path through the integration \cite{loffeld2013compar}. %and since the overhead of using \texttt{chop} does not allow us to use wall-clock time.
Since the sparse matrix-vector multiply is typically a memory bound computation, we set $a$ to be the ratio of double and single memory bandwidth and $b$ to be the ratio of double and half memory bandwidth. 
For typical hardware, like the NVIDIA A100, this simply yields $a = 2$ and $b = 4$.
This estimate may be conservative if using lower-precision moves the sparse matrix-vector multiply into a compute-bound regime (possible on some hardware, like the Cerebras Wafer Scale Engine \cite{jacquelin2022scalable}).
The notable result is that the mixed precision IOP Arnoldi makes all of the schemes more efficient in most regimes. 
The few exceptions are in the case of exprk4s6 when the error is around $10^{-10}$. 
In this case, the extra Krylov iterations induced by the lower precision introduce too much overhead for the use of low precision to provide a benefit.

\section{Conclusions}
\label{s:conclusion}

Modern computer hardware offers significantly increased low precision floating point performance in comparison to double precision.
We have developed two approaches to leveraging low precision in exponential time integration methods.

With a minor modification to the exponential Rosenbrock--Euler method, our reformulated version \eqref{eq:mixed_precision_exp_Euler} of the method attains improved resilience to inexact $\varphi$-functions. This enables utilizing cheaper, low precision arithmetic or looser tolerances in the most expensive part of the integrator. The reformulated exponential Rosenbrock--Euler method \eqref{eq:mixed_precision_exp_Euler} is particularly effective at maintaining convergence when combining single precision $\varphi$-functions with double precision for the remaining computations. Half precision $\varphi$-functions present many challenges with avoiding underflow and overflow, particularly for stiff problems. Nevertheless, improved accuracy is still achievable with our reformulated scheme \eqref{eq:mixed_precision_exp_Euler}. While we focused on the exponential Rosenbrock--Euler method, the reformulation idea could be generalized to other exponential methods and will the subject of future investigations.

Our mixed precision IOP Arnoldi algorithm incorporated into \KIOPS{}, or similar Krylov approximation methods, can readily be utilized within higher order methods as demonstrated in experiments with the advection-diffusion-reaction PDE.
This algorithm enables exponential methods to compute the $\varphi$-function products while leveraging low precision for the matrix-vector products and still recovering the required approximation accuracy.
The process of manually choosing $\mchop{1,2}$ for fixed matrices as in Section \ref{s:KIIOPS} is much more difficult in the context of ODE or spatially discretized PDE systems like the the advection-diffusion-reaction problem. 
This is primarily due to the dynamical nature of the problem changing the optimal values.
As such, to make the mixed precision IOP Arnoldi more practical, an adaptive approach to selecting the precision for the matrix-vector products is needed.
This is another topic we will explore in the future.

\section{Acknowledgements}

We would like to thank Valentin Dallerit for enlightening discussions and insight into the \KIOPS{} algorithm and software implementation. 
We are also grateful for the thoughtful feedback and insight of Professor Daniel Reynolds at Southern Methodist University.

% \appendix

\bibliographystyle{siam}
\bibliography{paper}

\end{document}